\documentstyle[12pt]{article}

\newtheorem{thm}{Theorem}
\newtheorem{lem}[thm]{Lemma}
\newtheorem{cor}[thm]{Corollary}
\newtheorem{pr}[thm]{Proposition}
\newtheorem{ex}[thm]{Example}
\newenvironment{pf}{\noindent {\em Proof:}}{$\Box$\\}

\newcommand{\N}{\mbox{\hskip.1em N \hskip -1.25em \relax I \hskip .1em}}
\newcommand{\R}{\mbox{\hskip.1em R \hskip -1.25em \relax I \hskip .1em}}
\newcommand{\J}{$J(e_i)$}
\newcommand{\Jd}{$J(e_i)'$}
\newcommand{\Jdd}{J(e_i)'}
\newcommand{\sumA}{\sum_{i\in A}}
\newcommand{\sumAc}{\sum_{i\notin A}}
\newcommand{\rnorm}{\biggr\|}
\newcommand{\lnorm}{\biggl\|}
\newcommand{\ei}{$(e_i)$}
\newcommand{\ui}{$(u_i)$}
\newcommand{\uid}{$(u'_i)$}
\newcommand{\eid}{$(e'_i)$}
\newcommand{\llang}{\biggl\langle}
\newcommand{\lrang}{\biggr\rangle}
\newcommand{\lv}{\biggl|}
\newcommand{\rv}{\biggr|}

\begin{document}

\begin{center}
  {\Large\bf  Banach Spaces with Property (w)}\vspace{3mm}\\
  {\large\sc Denny H. Leung}
\end{center}

\vspace{1mm}

\begin{abstract}
  A Banach space $E$\/ is said to have Property (w) if every (bounded
linear) operator from $E$\/ into $E'$\/ is weakly compact.  We give some 
interesting
examples of James type Banach spaces with Property (w). We also consider the 
passing
of Property (w) from $E$\/ to $C(K,E)$.
%abstract
\end{abstract}

\begin{figure}[b]
  \rule{3in}{.005in}
\end{figure}

\baselineskip 4ex

\section{Introduction}
A Banach space $E$\/ is said to have Property (w) if every operator
from $E$\/ into $E'$\/ is weakly compact.  This property was
introduced by E.\ and P.\ Saab in \cite{SS}.  They observe that for
Banach lattices, Property (w) is equivalent to Property (V*), which in
turn is equivalent to the Banach lattice having a weakly sequentially
complete dual.  Thus the following question was raised \cite{SS}:\\

\noindent{\em Question}: Does every Banach space with Property (w) have a
weakly sequentially complete dual, or even Property (V*)? \\

In this paper, we give two examples, both of which answer the question
in the negative.  Both examples are James type
spaces considered in \cite{BHO}.  They both possess properties
stronger than Property (w).  The first example has the property that
every operator from the space into the dual is compact.  In the second
example, both the space and its dual have Property (w).  In the last
section, 
we consider if Property (w) passes from a Banach space $E$\/ to
$C(K,E)$\/.  This was also dealt with in \cite{SS}.\\
\indent We use standard Banach space terminology as may be found in \cite{LT}.
For a Banach space $E$, $E'$ denotes its dual, and $U_E$ its closed
unit ball.  If $F$ is also a Banach space, then we let $L(E,F)$
(respectively $K(E,F)$) denote the space of all bounded linear
operators (respectively all compact operators) from $E$ into $F$.  The
norm in $\ell^p$ is denoted by $\|\cdot\|_p$.  Let us also establish
some terminology about sequences.  If $(e_i)$ is a sequence in a
Banach space, we use $[e_i]$ to denote the closed linear span of
$(e_i)$. The sequence is {\em semi-normalized\/} if $0 < \inf\|e_i\| \leq
\sup\|e_i\| < \infty$.  If $(f_i)$ is another sequence in a possibly
different Banach space, we say that $(e_i)$ {\em dominates} $(f_i)$ if there
is a constant $C$ such that $\|\sum a_if_i\| \leq C\|\sum a_ie_i\|$ for all
finitely non-zero real sequences $(a_i)$.  Two sequences are {\em
equivalent\/} if each dominates the other.  Finally, we use the symbol
$\preceq$ (respectively $\succeq$) to indicate $\leq$ (respectively
$\geq$) up to a fixed constant.  The symbol $\approx$ stands for
``$\preceq$ and $\succeq$''.

\section{James type constructions}

We first recall the construction of James type spaces as in
\cite{BHO}.  Let $(e_i)$ be a normalized basis of a Banach space $E$.
For $(a_i) \in c_{00}$, the space of all finitely non-zero real
sequences, let
\[ \biggl|\biggl|\biggl|\sum a_iu_i\biggr|\biggr|\biggr| =
\sup\biggl\{\lnorm\sum^k_{i=1}\biggl(\sum^{q(i)}_{j=p(i)}a_j\biggr)e_{p(i)}
\rnorm
: k \in \N\, ,
1\leq p(1) \leq q(1) < \ldots < p(k) \leq q(k)\biggr\}. \]
The completion of the linear span of the sequence $(u_i)$ is denoted
by $J(e_i)$.  Since we will only consider unconditional $(e_i)$, we
use the equivalent norm
\[ \lnorm\sum a_iu_i\rnorm =
\sup\biggl\{\lnorm\sum^k_{i=1}\biggl(\sum^{p(i+1)-1}_{j=p(i)}a_j
\biggr)e_{p(i)}\rnorm :
 k \in \N\, ,
1= p(1) < p(2) < \ldots < p(k+1)\biggr\}. \]
As in \cite{BHO}, $\sum^k_{i=1}(\sum^{p(i+1)-1}_{j=p(i)}a_j)e_{p(i)}$
is called a {\em representative} of $\sum a_iu_i$ in $E$.  The
biorthogonal sequences of \ei\ and $(u_i)$ are denoted by $(e'_i)$ and
$(u'_i)$ respectively.  The basis projections with respect to the
basis \ui\ are denoted by $(P_n)^\infty_{n=0}$\ $(P_0 = 0)$.  The
functional $S$ defined by $S(\sum a_iu_i) = \sum a_i$ is bounded,
hence $S \in \Jdd$.  The
following lemma is useful for computing the norms of certain
vectors in \J.  Recall that a basic sequence $(x_i)$ is {\em right
dominant} \cite{BHO} if it is unconditional and whenever $1 \leq m(1)
\leq n(1) < \ldots < m(i) \leq n(i) < \ldots$, we have
\[ \lnorm\sum a_{n(i)}x_{m(i)}\rnorm \preceq
\lnorm\sum a_{n(i)}x_{n(i)}\rnorm. \]
\vspace*{0ex}
\begin{lem}\label{norm}
Let $(e_i)$ be a right dominant basis of a Banach space.  There
exists a constant $C$ such that whenever $(v_i)^l_{i=1}$ is a block
basis of $(u_i)$ satisfying
\begin{equation}\label{vi}
v_i = \sum^{n(i+1)-1}_{j=n(i)}a_ju_j\, ,\hspace{1cm}
\sum^{n(i+1)-1}_{j=n(i)}a_j
= 0,\hspace{1cm} 1 \leq i \leq l ,
\end{equation}
there is a block basis $(w_i)^l_{i=1}$ of $(e_i)$ such that $||w_i||
\leq ||v_i||, 1 \leq i \leq l$, and $||\sum^l_{i=1}v_i|| \leq
C||\sum^l_{i=1}w_i||$.
\end{lem}

\begin{pf}
Let $(v_i)$ be as given.  Choose $p(1) < p(2) < \ldots < p(k+1)$ such
that 
\[ \lnorm\sum^l_{i=1}v_i\rnorm \preceq
\lnorm\sum^k_{i=1}\biggl(\sum^{p(i+1)-1}_{j=p(i)}a_j\biggr)e_{p(i)}\rnorm. \]
Let $(q(i))$ be a finite strictly increasing sequence such that 
\[ \{q(i)\} = \{p(i): 1 \leq i \leq k+1\} \cup \{n(i): 1 \leq i \leq
l+1\}. \]
Let $A = \{i:$\/ there exists some $m$\/ with $p(i) < n(m) <
p(i+1)\}$.  For
all $i \in A$, let $m_i = \min\{m:n(m)>p(i)\}$ and $r_i =
\max\{m:n(m)<p(i+1)\}$.  By (\ref{vi}), for $i \in A$,
\[ \sum^{p(i+1)-1}_{j=p(i)}a_j = \sum^{n(m_i)-1}_{j=p(i)}a_j +
\sum^{p(i+1)-1}_{j=n(r_i)}a_j \equiv b_i + c_i .\]
Therefore,
\begin{eqnarray*}
\lnorm\sum_{i\in A}(\sum^{p(i+1)-1}_{j=p(i)}a_j)e_{p(i)}\rnorm & \leq
& \lnorm \sumA
b_ie_{p(i)}\rnorm + \lnorm\sumA c_ie_{p(i)}\rnorm \\
& \preceq & \lnorm\sumA b_ie_{p(i)}\rnorm + \lnorm\sumA
c_ie_{n(r_i)}\rnorm \\& & 
\hspace{2cm} \mbox{\rm since $(e_i)$ is right dominant}\\
& \preceq & \lnorm\sum\biggl(\sum^{q(i+1)-1}_{j=q(i)}a_j\biggr)e_{q(i)}\rnorm.
\end{eqnarray*}
Also, it is clear that
\[\lnorm\sumAc\biggl(\sum^{p(i+1)-1}_{j=p(i)}a_j\biggr)e_{p(i)}\rnorm
 \preceq \lnorm\sum\biggl(\sum^{q(i+1)-1}_{j=q(i)}a_j\biggr)e_{q(i)}\rnorm. \]
Hence 
\[ \lnorm\sum^l_{i=1} v_i\rnorm \preceq
\lnorm\sum\biggl(\sum^{q(i+1)-1}_{j=q(i)}a_j\biggr)e_{q(i)}\rnorm. \]
Now let $B_m = \{i:n(m) \leq q(i) < n(m+1) \}$.  Then
\begin{eqnarray*}
\lnorm\sum^l_{i=1}v_i\rnorm & \preceq & \lnorm\sum_m\biggl(\sum_{i\in
B_m}(\sum^{q(i+1)-1}_{j=q(i)}a_j)e_{q(i)}\biggr)\rnorm \\
& \equiv & \lnorm\sum w_m\rnorm.
\end{eqnarray*}
Note that $w_m$ is a representative of $v_m$ in $E$.  Hence $||w_m||
\leq ||v_m||$.  Clearly, $(w_m)$ is a block basis of $(e_i)$.
\end{pf}

\begin{lem}\label{perturb}
Let \ei\ be a shrinking normalized unconditional basic sequence.
Assume that $L(E,\Jdd) \neq K(E,\Jdd)$ for some
subspace $E$ of \J. Then there are semi-normalized block bases
$(x_i)$ and $(x'_i)$ of \ui\ and $(u'_i)$ respectively such that
$(x_i)$ dominates $(x'_i)$ and $Sx_i = 0$ for all $i$.
\end{lem}

\begin{pf}
Since \ei\ is shrinking, by \cite{BHO}, Theorems 2.2
and 4.1, \Jd $= [\{S\} \cup \{u'_i\}^\infty_{i=1}]$.  
Let $T:E\rightarrow\Jdd$ be non-compact. Define the projection
$P:\Jdd\rightarrow[u'_i]$ by $P(aS+\sum a_iu'_i) = \sum
a_iu'_i$.  Then $(1-P)T$ has rank $1$.  Thus $PT$ is non-compact.
Replacing $T$ by $PT$, we may assume without loss of generality that
range\,$T \subseteq [u'_i]$.  Choose a bounded sequence $(y_i)$ in $E$
so that $\inf_{i,j}\|Ty_i-Ty_j\| > 0$.  By \cite{BHO}, Theorem 4.1,
$\ell^1$ does not embed into \J.  Hence we may assume that $(y_i)$ is
weakly Cauchy.  Thus $(y_{2i-1}-y_{2i})$ is weakly null and
semi-normalized; hence by Proposition 1.a.12 of \cite{LT}, we may
assume that it is equivalent to a semi-normalized block basis $(x_i)$
of \ui.  Since $(x_i)$ is weakly null, $Sx_i \to 0$.  By using
a subsequence, we may further assume that $Sx_i = 0$ for all $i$.
Similarly, $(T(y_{2i-1}-y_{2i})) \subseteq [u'_i]$ is semi-normalized
and weakly Cauchy.  Using the same argument, we may assume that it is
equivalent to a semi-normalized block basis  $(x'_i)$ of $(u'_i)$.
Finally,
\begin{eqnarray*}
\lnorm\sum a_ix'_i\rnorm & \approx & \lnorm\sum a_iT(y_{2i-1}-y_{2i})\rnorm \\
& \leq & \|T\| \, \lnorm\sum a_i(y_{2i-1}-y_{2i})\rnorm \\
& \approx & \lnorm\sum a_ix_i\rnorm,
\end{eqnarray*}
as required.
\end{pf} 
 
\begin{thm}\label{compact}
Let $(e_i)$ be a right dominant, normalized basic sequence which
dominates all of its normalized block bases.  Then the following are
equivalent.\\
{\rm (a)}  For every subspace $E$ of \J, $L(E,\Jdd) = K(E,\Jdd)$,\\
{\rm (b)}  The sequence $(e_i)$ does not dominate $(e'_i)$.
\end{thm}

\begin{pf}
Since $(e_i)$ is unconditional, if it is not shrinking, it has a
normalized block basis equivalent to the $\ell^1$ basis.  Hence \ei\
dominates the $\ell^1$ basis.  Thus \ei\ is equivalent to the $\ell^1$
basis.  Both (a) and (b) fail in this case, so they are equivalent.
Now assume that \ei\ is shrinking.\\

\noindent (a)$\Rightarrow$(b). 
 Assume \ei\ dominates \eid.  If $||\sum b_iu_i||
\leq 1$, then $||\sum b_ie_i|| \leq 1$.  Hence
\begin{eqnarray*}
\lv\llang\sum b_iu_i, \sum a_iu'_{2i}\lrang\rv & = & \lv\sum a_ib_{2i}\rv \\
& = & \lv\llang\sum b_ie_i\, , \sum a_i e'_{2i}\lrang\rv \\
& \leq & \lnorm\sum a_ie'_{2i}\rnorm.
\end{eqnarray*}
Therefore, $(u'_{2i})$ is dominated by $(e'_{2i})$.  By \cite{BHO},
Proposition 2.3, $(u_{2i-1}-u_{2i})$ is equivalent to $(e_{2i})$.
Hence, if $E = [u_{2i-1}-u_{2i}]$, then the map $T:E\rightarrow\Jdd$ defined
by  $T(\sum a_i(u_{2i-1}-u_{2i})) = \sum a_iu'_{2i}$ is bounded.
Clearly, $T$ is not compact.  \\

\noindent (b)$\Rightarrow $(a).  Suppose $L(E,\Jdd) \neq K(E,\Jdd)$ for some
subspace $E$ of \J.
By Lemma \ref{perturb}, there are semi-normalized block bases
$(x_i)$ and $(x'_i)$ of \ui\ and $(u'_i)$ respectively such that
$(x_i)$ dominates $(x'_i)$ and $Sx_i = 0$.  Let $1\leq n_1 \leq m_1 <
n_2 \leq m_2 < \ldots$ be such that $x'_i \in [u'_n]^{m_i}_{n=n_i}$.
By using a subsequence, assume that
\begin{equation}\label{sep}
m_i - n_i +1 \leq n_{i+1} - m_i -1 \hspace{1em}\mbox{for
all}\hspace{1em} i \geq 1.
\end{equation}
Let $R:$\J$\to$\J\ denote the right shift operator $R(\sum a_iu_i) =
\sum a_iu_{i+1}$.  Since \ei\ dominates all of its normalized block
bases, $(R^k)$ is uniformly bounded.  Choose a normalized block basis
$(y_i)$ of \ui\ such that $y_i \in [u_n]^{m_i}_{n=n_i}$ and $\langle
y_i,x'_i\rangle = \|x'_i\|$ for all $i$.  Now let $z_i = y_i -
R^{m_i+1-n_i}y_i$ for all $i \geq 1$.  By (\ref{sep}), $(z_i)$ is a
semi-normalized block basis of \ui.  Also $\langle
z_i,x'_i\rangle = \|x'_i\|$ and $Sz_i = 0$ for all $i$. Fix $(a_i) \in
c_{00}$.  By Lemma
\ref{norm}, there is a block basis $(w_i)$ of \ei\ such that $\|w_i\|
\leq \|a_iz_i\|$ for all $i$, and
\begin{eqnarray*}
\lnorm\sum a_iz_i\rnorm & \preceq & \lnorm\sum w_i\rnorm \\
& = & \lnorm\sum\|w_i\|\frac{w_i}{\|w_i\|}\rnorm \\
& \preceq & \lnorm \sum\|w_i\|e_i\rnorm \hspace{2em}\mbox{since \ei\ dominates }
\biggl(\frac{w_i}{\|w_i\|}\biggr) \\
& \preceq & \lnorm\sum a_ie_i\rnorm .
\end{eqnarray*}
Computing the norm of a vector of the form $\sum b_ix'_i$ on some $\sum
a_iz_i$, we find that $(x'_i)$ dominates
\eid.  Since $Sx_i = 0$ for all $i$.  The computation above can also
be applied to $(x_i)$.  Hence $(x_i)$ is
dominated by \ei.  Since $(x_i)$ dominates $(x'_i)$ by choice of the
sequences, we see that \ei\ dominates \eid.
\end{pf}

\begin{cor}
Let \ei\ be the unit vector basis of $\ell^p, 2 < p < \infty$.  Then
\J\ has Property~{\em(w)} but not a weakly sequentially complete dual.
\end{cor}

\begin{pf}
This follows immediately from Theorem \ref{compact} and the fact that
\J\ is quasi-reflexive of order 1 by \cite{BHO}, Theorem 4.1.
\end{pf}

\noindent{\em Remark}.  Theorem \ref{compact} fails without the
assumption that \ei\ dominates all of its normalized block bases.  In
fact, if \ei\ is subsymmetric, then by \cite{BHO}, Propostion 2.3,
\ei\ is equivalent to $(u_{2i-1}-u_{2i})$.  Similarly, \eid\ is
equivalent to $(u'_{2i})$ in this case.  Now if we let \ei\ be the
unit vector basis of the Lorentz space $d(w,2)$ \cite{LT}, then \ei\
is symmetric and shrinking, and does not dominate \eid.  However,
$\ell^2$ embeds into both $d(w,2)$ and its dual, and hence into both
\J\ and \Jd\ by the observation made above.  Hence condition (a) of
Theorem \ref{compact} fails.

\section{A non-reflexive space whose dual and itself have Property (w)}

In this section, we give an example of a non-reflexive Banach space
$E$ so that both $E$ and $E'$ have Property (w).  In fact, neither $E$
nor any of its higher duals is weakly sequentially complete.  The
example will again be a \J\ space with a suitably chosen $(e_i)$.\\
\indent If \ei\ is a normalized basis of a reflexive Banach space $E$, then
\J\ is quasi-reflexive of order 1 by \cite{BHO}, Theorem 4.1.  Thus,
if we define the functional $L$ on \Jd\ by $L(aS+\sum a_iu'_i) = a$,
then $J(e_i)'' = [\{L\}\cup\{u_i\}^\infty_{i=1}]$.  Using this
observation, the following Proposition can be obtained by straight
forward perturbation arguments.

\begin{pr}\label{wperturb}
Let \ei\ be a normalized unconditional basis of a reflexive Banach
space $E$.\\
\noindent{\rm (a)} Let $(y_n)$ be a bounded sequence in \J\ with no weakly
convergent subsequence.  Then there exist a subsequence $(y_{n_i})$,
an element $y_0$ of \J, a
block basis $(z_i)$ of $(u_n)$, and $a \neq 0$ such that\\ 
\indent{\rm (i)} $(y_{n_i}-y_0) \approx (z_i)$, and\\
\indent{\rm (ii)} $Sz_i = a$ for all $i$.\\
\noindent{\rm (b)} Let $(y'_n)$ be a bounded sequence in $[u'_n]$ with
no weakly convergent subsequence. Then there are a subsequence
$(y'_{n_i})$, a vector $y'_0 \in [u'_n]$, $0 = k_0 < k_1 < \ldots,
(z'_i) \subseteq [u'_n]$, and $b \neq 0$ such that\\
\indent{\rm (i)} For all $i$, $z'_i \in [u'_n]^{k_i}_{n=k_{i-1}+1}$, and\\
\indent{\rm (ii)} The sequence $(y'_{n_i} - y'_0 - bP'_{k_{i-1}}S - z'_i)$
is norm null. 
\end{pr}

\begin{pr}\label{W}
Let \ei\ be a subsymmetric normalized basis of a reflexive Banach
space $E$, and let $1<p<\infty$. Assume that\\
\noindent{\rm (i)} $E$ satifies an upper $p$-estimate,\\
\noindent{\rm (ii)} \eid\ does not dominate \ei, and\\
\noindent{\rm (iii)} The $\ell^p$ basis does not dominate \eid.\\
\noindent Then both \J\ and \Jd\ have Property {\em (w)}, but neither has a
weakly sequentially complete dual.
\end{pr}

\begin{pf}
Since \J\ is quasi-reflexive of order 1, neither \J\ nor \Jd\ has a
weakly sequentially complete dual.  If $T: $\J\ $\to\Jdd$\ is not weakly
compact, then we may assume that range\,$T \subseteq [u'_i]$ as in the
proof of Lemma \ref{perturb}.  Now there is a bounded sequence $(y_n)$
in \J\ such that $(Ty_n)$ has no weakly convergent subsequence and 
$\inf_{i,j}\|y_i-y_j\| > 0 $.  Apply Proposition
\ref{wperturb} to $(y_n)$ and $(y'_n) \equiv (Ty_n)$ to yield the various
objects identified there. For all $i$, let $x_i = y_{n_{4i-1}}-y_{n_{4i-3}}$.
Then $(x_i)$ is semi-normalized and $(x_i) \approx (z_{4i-1}-z_{4i-3})$ 
(where $(z_i)$ is as given by Proposition \ref{wperturb}). 
Also $\|Tx_i-x'_i\|\to 0$, where
\[ x'_i = b(P'_{k_{4i-1}}-P'_{k_{4i-3}})S+z'_{4i-1}-z'_{4i-3}.\]
Note that $(x'_i)$ is a semi-normalized block basis of \uid\ such that
$x'_i \in [u'_n]^{k_{4i-1}}_{n=k_{4i-4}+1}$ and
$\langle u_j, x'_i\rangle = b$ if $k_{4i-3}<j\leq k_{4i-2}$.  Without
loss of generality, assume $(Tx_i)\approx (x'_i)$.  For all $i$, let
$k_{2i-1}<j_i\leq k_{2i}$.  By the subsymmetry of \ei, it is easy to
see that $\|\sum b_i(u_{j_{2i-1}}-u_{j_{2i}})\| \approx \|\sum b_ie_i\|$.
Computing the norm of a vector of the form $\sum a_ix'_i$ on $\sum
b_i(u_{j_{2i-1}}-u_{j_{2i}})$, we see that $(x'_i)$ dominates \eid.  On
the other hand, since $S(z_{4i-1}-z_{4i-3}) = 0$ for all $i$. Lemma
\ref{norm} implies that $(x_i) \approx (z_{4i-1}-z_{4i-3})$ is
dominated by some semi-normalized block basis of \ei.  But since $E$
satisfies an upper $p$-estimate, $(x_i)$ is dominated by the $\ell^p$
basis.  Thus
\begin{eqnarray*}
\lnorm \sum a_ie'_i\rnorm &\preceq & \lnorm\sum a_ix'_i\rnorm \\
& \approx & \lnorm\sum a_iTx_i\rnorm \\
& \preceq & \lnorm\sum a_ix_i\rnorm \\
& \preceq & \biggl(\sum|a_i|^p\biggr)^{1/p},
\end{eqnarray*}
a contradiction.  Hence \J\ has Property (w).\\

Since \J\ is quasi-reflexive, if \Jd\ fails Property (w), there is an
operator $T~:~[u'_n]~\to$~\J\ which is not weakly compact.  Choose a bounded
sequence $(y'_n)$ in $[u'_n]$ so that $(y_n)\equiv (Ty'_n)$ has no weakly
convergent subsequence.  As before, apply Lemma \ref{wperturb}.
In the present situation, we may assume $z'_i = 0$ as $(z'_i)$ is a
weakly null sequence.  Arguing as before, we find that the sequence
$(z_{2i}-z_{2i-1})$ is dominated by $((P'_{k_{2i}}-P'_{k_{2i-1}})S)$.
Since $Sz_i = a \neq 0$ for all $i$, and because \ei\ is
subsymmetric, it follows that $(z_{2i}-z_{2i-1})$ dominates \ei.  On
the other hand, if $x\in$ \J, $\|x\|\leq 1$, then
$\|\sum\langle(P_{k_{2i}}-P_{k_{2i-1}})x,S\rangle e_i\|\preceq 1$.
Therefore,
\begin{eqnarray*}
\lv\llang x, \sum a_i(P'_{k_{2i}}-P'_{k_{2i-1}})S\lrang\rv & = &
\lv\sum a_i\llang(P_{k_{2i}}-P_{k_{2i-1}})x,S\lrang\rv \\
& \preceq & \lnorm\sum a_ie'_i\rnorm.
\end{eqnarray*}
Hence \eid\ dominates $((P'_{k_{2i}}-P'_{k_{2i-1}})S)$.  Consequently,
\eid\ dominates \ei.  This contradicts assumption (ii).  Hence \Jd\
also has Property (w).
\end{pf}  
 
We now construct a sequence \ei\ satisfying the conditions in
Propostion \ref{W}. For a number $p \in (1,\infty)$, let $p' =
p/(p-1)$.  Now fix $p$ and $r$ such that $1 < p < 2 < r < \infty$ and
$r' < p$.  Let $(k_n) \subseteq \N\ $ be such that
\begin{equation}\label{sum}
 1 + 2 \sum^{n-1}_{i=1}i^{r/2}k^{r/2}_i \hspace{1em}\leq\hspace{1em}
k^{1-p/2}_n \bigg/
n^{p/2}\hspace{2em}\mbox{and} 
\end{equation}
\begin{equation}\label{frac}
\biggl(\frac{n+1}{n}\biggr)^{r/2}\biggl(\frac{k_n}{k_{n+1}}
\biggr)^{(1/p-1/p')r/2}\hspace{1em}
\leq\hspace{1em} \frac{1}{2}\ .
\end{equation}
Finally, let $\alpha_n = \sqrt{n}k_n^{(1/p'-1/p)/2}$ for all $n \geq
1$.\\

\begin{lem}\label{calc}
For all $l, j \in$ {\rm \N\ }  with 
\begin{equation}\label{j}
 2 \sum^{l-1}_{i=1}i^{r/2}k^{r/2}_i \hspace{1em}\leq\hspace{1em} j
\hspace{1em}\leq\hspace{1em} k^{1-p/2}_l \bigg/ l^{p/2}, 
\end{equation}
and for all $(x_i)^l_{i=1}$ with $0 \leq x_i \leq k_i, 1 \leq i < l$,
$0 \leq x_l \leq  j - \sum^{l-1}_{i=1} x_i $, we have
\[0 \hspace{1em}\leq\hspace{1em} \sum^l_{i=1} \alpha^r_i x^{r/p}_i
\hspace{1em}\leq\hspace{1em} \frac{j}{2} + 2. \]
\end{lem}

\begin{pf}
Note that $j - \sum^{l-1}_{i=1} k_i \geq 0$ by the choice of $j$,
since $r > 2$.  Clearly, there is no loss of generality in assuming
that $x_l = j - \sum^{l-1}_{i=1} x_i$.
Let $f:\prod^{l-1}_{i=1}[0,k_i]\to \R\ $be defined by 
\[ f(x_1,\ldots,x_{l-1}) = \sum^l_{i=1}\alpha^r_ix^{r/p}_{i} .\]
Clearly, $f \geq 0$ on $\prod^{l-1}_{i=1}[0,k_i]$.  Suppose
that $f$ attains its maximum at $(j_1,\ldots,j_{l-1}) \in
\prod^{l-1}_{i=1}[0,k_i]$.  Let $A = \{1 \leq i < l : j_i = 0\}$ and $B
= \{1 \leq i < l : j_i = k_i\}$.  Then for $i \notin A \cup B$, 
\begin{eqnarray*}
\frac{\partial f}{\partial x_i}(j_1,\ldots,j_{l-1}) & = & 0 \\ 
\Rightarrow \hspace{5em} \alpha^r_ij^{r/p}_i & = & \alpha^r_lj^{r/p-1}_lj_i .
\end{eqnarray*}
Hence, 
\[ f(j_1,\ldots,j_{i-1}) = \sum_{i \in B}\alpha^r_ik^{r/p}_i +
\sum_{\begin{array}{c}\\[-3.6ex]{\scriptstyle 1 \leq i < l} \\ 
\\[-3.3ex]{\scriptstyle i\notin A \cup B}\end{array}}\!\! 
\alpha^r_lj^{r/p-1}_lj_i + \alpha^r_lj^{r/p}_l , \]
where $j_l = j - \sum^{l-1}_{i=1}j_i$.  Now $\alpha^r_ik^{r/p}_i =
i^{r/2}k^{r/2}_i$ .  Hence
\begin{eqnarray*}
f(j_1,\ldots,j_{l-1}) & = & \sum_{i \in B}i^{r/2}k^{r/2}_i +
\alpha^r_lj^{r/p}_l\biggl(1 + 
\sum_{\begin{array}{c}\\[-3.6ex]{\scriptstyle 1 \leq i < l} \\
\\[-3.3ex]{\scriptstyle i \notin A \cup B}
\end{array}}\!\! j_i\bigg/j_l \biggr) \\
& \leq & \frac{j}{2} +
\alpha^r_lj^{r/p}_l\biggl(1 + \sum^{l-1}_{i=1}j_i\bigg/j_l
\biggr) \\
& \leq & \frac{j}{2} +
\alpha^r_lj^{r/p}_l\biggl(1 +
\sum^{l-1}_{i=1}k_i\bigg/\biggl(j-\sum^{l-1}_{i=1}k_i\biggr)\biggr) \\
& \leq & \frac{j}{2} +
2\alpha^r_lj^{r/p}_l\hspace{2em}\mbox{by choice of $j$} \\
& \leq & \frac{j}{2} + 2\alpha^r_lj^{r/p} \\
& \leq & \frac{j}{2} +
2\biggl(\alpha_lk^{(1/p-1/2)}_l\bigg/l^{1/2}\biggr)^r \hspace{2em}
\mbox{by choice of $j$} \\
& = & \frac{j}{2} + 2 .
\end{eqnarray*}
\end{pf}

Let $(t_i)$ be the normalized basic sequence so that for all
$(a^i_j)^{k_i}_{j=1}\mbox{}^\infty_{i=1} \in c_{00}$,
\[ \lnorm\sum^\infty_{i=1}\biggl(\sum^{k_i}_{j=1}
a^i_jt_{k_1+\ldots+k_{i-1}+j}\biggr)\rnorm
= \|(a^i_j)\|_r \vee
\lnorm\biggl(\alpha_i\biggl(\sum^{k_i}_{j=1}|a^i_j|^p\biggr)^{1/p}
\biggr)^\infty_{i=1}\rnorm_r.~~~(k_0
= 0) \]
\vspace*{0ex}
\begin{ex}
Let $E$ be the Banach space with a basis \ei\ so that
\[ \lnorm\sum b_ie_i\rnorm = \sup\biggl\{\lnorm\sum b_it_{\sigma(i)}\rnorm
: \mbox{$\sigma$ is a permutation of\ ${\rm \N\ }$}\biggr\} \]
Then $E$ satisfies the conditions of Proposition\ {\rm \ref{W}}.  Hence \J\
and \Jd\ have Property\ {\rm (w)}, but neither has a weakly sequentially
complete dual.
\end{ex}

\begin{pf}
It is clear that \ei\ is a symmetric basis and $E$ satisfies an upper
$p$-estimate.  To show that $E$ is reflexive, it suffices to show that
$c_0$ does not lattice embed into $E$.  This will follow if we show
that $c_0$ does not lattice embed into $[t_n]$.  Let $(x_i)$ be a disjoint
normalized sequence in $[t_n]$.  If $\inf\|x_i\|_r = \epsilon > 0$, then
\[ \lnorm\sum^k_{i=1}x_i\rnorm \geq \lnorm\sum^k_{i=1}x_i\rnorm_r \geq
\biggl(\sum^k_{i=1}\|x_i\|^r_r\biggr)^{1/r} \geq \epsilon
k^{1/r} \]
for all $k$.  Hence $(x_i)$ is not equivalent to the $c_0$ basis.  On
the other hand, if $\inf\|x_i\|_r = 0$, we may assume that $\|x_i\|_r\to
0$.  For each $i$, we write $x_i = \sum_j x_i(j)$, with $x_i(j) \in
[t_n]^{k_1+\ldots+k_j}_{n=k_1+\ldots+k_{j-1}+1}$. Since $\|x_i\|_r\to
0$, $\lim_i x_i(j) = 0$ for all $j$.  Thus, by perturbation and dropping to
a subsequence, we may assume that each $x_i$ has the form
$\sum^{j_i}_{j=j_{i-1}+1} x_i(j)$, where $(j_i)$ is strictly
increasing.  But then since $(x_i)$ is normalized and these factors
add according to the $\ell^r$ norm, we see that $(x_i)$ is not equivalent
to the $c_0$ basis.  This shows that $E$ is reflexive.\\
\indent For all $l$, choose $j$ as in Equation (\ref{j}).  We estimate
the norm of
$\sum^j_{i=1} e_i$.  It is easy to see that there exists
$(j_i)^\infty_{i=1} \subseteq \N\ \cup\{0\}$ such that $\sum j_i = j$
and 
\[ \lnorm\sum^j_{i=1}e_i\rnorm =  j^{1/r} \vee
\lnorm(\alpha_ij^{1/p}_i)^\infty_{i=1}\rnorm_r .\]
Now 
\[ \sum^l_{i=1}\alpha^r_ij^{r/p}_i \leq \frac{j}{2} + 2 \]
by Lemma \ref{calc}.  For $i > l$, 
\begin{eqnarray*}
\alpha_ij^{1/p}_i & \leq & \alpha_ij^{1/p} \\
& \leq & \alpha_i(k^{1-p/2}_l/l^{p/2})^{1/p} \\
& = & \sqrt{\frac{i}{l}}\biggl(\frac{k_l}{k_i}\biggr)^{(1/p-1/p')/2} \\
& = &
\biggl(\prod^{l+1}_{m=i}\frac{m}{m-1}\biggr)^{1/2} 
\biggl(\prod^{l+1}_{m=i}\frac{k_{m-1}}{k_{m}}\biggr)^{(1/p-1/p')/2}\\
& \leq & 2^{-(i-l)/r},
\end{eqnarray*}
by Equation (\ref{frac}).  It follows easily that
$\|(\alpha_ij^{1/p}_i)^\infty_{i=1}\|_r \preceq j^{1/r}$.  Hence
$\|\sum^j_{i=1}e_i\| \preceq j^{1/r}$.  This implies that
\[ \lnorm\sum^j_{i=1}e'_i\rnorm \succeq j^{1/r'}. \]
Since $r' < p$, the $\ell^p$ basis does not dominate $(e'_i)$.\\
\indent Finally, for all $n$,
\[ \lnorm\sum^{k_n}_{i=1}e_i\rnorm \geq \alpha_nk^{1/p}_n =
\sqrt{nk_n} .\]
Also,
\begin{eqnarray*}
\lnorm\sum^{k_n}_{i=1}b_ie_i\rnorm & \leq & 1\\
\Rightarrow\hspace{3.3em} \alpha_n\biggl(\sum^{k_n}_{i=1}|b_i|^p\biggr)^{1/p} &
\leq & 1 \\
\Rightarrow\hspace{3.1em}\llang\sum^{k_n}_{i=1}b_ie_i ,
\sum^{k_n}_{i=1}e'_i\lrang &
\leq & \|(b_i)\|_p\, k^{1/p'}_n \\
& \leq & \frac{k^{1/p'}_n}{\alpha_n} \\
& = & \sqrt{\frac{k_n}{n}} .
\end{eqnarray*}
Hence $\|\sum^{k_n}_{i=1}e'_i\| \leq (k_n/n)^{1/2}$.  Therefore, \eid\
does not dominate \ei.
\end{pf}

\section{Property (w) in $C(K,E)$}

In this section, we consider the question of whether the Property (w)
passes from a Banach space $E$ to $C(K,E)$, the space of all
continuous $E$-valued functions on a compact Hausdorff space $K$.  Let
$E, F$ be Banach spaces, and let $K$ be an arbitrary compact Hausdorff
space whose collection of Borel subsets is denoted by $\Sigma$. The
dual of $C(K,E)$ is isometric to the space $M(K,E')$ of all regular
$E'$-valued measures of bounded variation on $K$ \cite{DU}.  In case
$E = \R$\, , we use the notation $C(K)$ and $M(K)$ respectively. For
$\mu$ in $M(K,E')$, let $|\mu| \in M(K)$ denote its variation
(\cite{DU}, p. 2). 
Given $T \in L(C(K,E),F)$, it is well known that $T$ can be represented by a
vector measure $G : \Sigma \to L(E,F'')$ \cite{DU}.  In fact, $G$ is given by
\[ G(A)x = T''(\chi_A\otimes x) \]
for all $A \in \Sigma$ and $x \in E$, where $\chi_A$ is the
characteristic function of the set $A$.  For all $y' \in F'$, we
define $G_{y'} \in M(K,E')$ by 
\[ \langle x, G_{y'}(A)\rangle = \langle y', G(A)x\rangle \]
 for all $x \in E, A
\in \Sigma$. Then the {\em semivariation} of $G$ is given by
\[ \|G\|(A) = \sup \{|G_{y'}|(A) :  \|y'\| \leq 1\}. \]
The following result is well known \cite{BL}.

\begin{pr}\label{op}
Let $T:C(K,E)\to F$ be weakly compact.  Then its representing measure
$G$ satisfies\\
{\rm (a)} $G$ takes values in $L(E,F)$,\\
{\rm (b)} For every $A \in \Sigma$, $G(A)$ is weakly compact, and\\
{\rm (c)} The semivariation of $G$ is continous at $\emptyset$, {\em i.e.},
$\lim_n\|G\|(A_n) = 0$ for every sequence $(A_n)$ in $\Sigma$ which decreases
to $\emptyset$.
\end{pr}

Condition (c) is just the
uniform countable additivity of the set 
$\{|G_{y'}|: \|y'\| \leq 1\}$.
By Lemma VI.2.13 of \cite{DU}, this is 
equivalent to the fact that 
\[ \lim_n \sup_{\|y'\|\leq 1} |G_{y'}|(A_n)= 0 \]
whenever $(A_n)$ is a pairwise disjoint sequence in $\Sigma$.  For
the sake of brevity, we introduce the following {\em ad hoc}
terminology.\\ 

\noindent{\bf Definition}.  A pair $(K,E)$, where $K$ is a compact Hausdorff
space and $E$ is a Banach space, is called {\em simple} if for every
Banach space $F$, every operator  $T:C(K,E)\to F$ whose
representing measure $G$ satifies conditions (a)--(c) of Proposition
\ref{op} is weakly compact.

\begin{thm}
Suppose $(K,E)$ is simple.  Then $E$ has Property {\rm (w)}\, if and
only if $C(K,E)$ does.
\end{thm}

\begin{pf}
One direction is trivial.  Now assume that $E$ has Property (w).  Then
$E'$ cannot contain a copy of $c_0$ \cite{SS}.  We first prove the\\

\noindent{\bf Claim}.  Every (bounded linear) operator from $C(K,E)$
into $E'$ is weakly compact.\\

Let $T:C(K,E)\to E'$ be represented by the  measure $G$. 
For all $x \in E$, define
$T_x:C(K)\to E'$ by $T_xf = T(f\otimes x)$ for all $f \in C(K)$.  Since
$E'$ does not contain a copy of $c_0$, $T_x$ is weakly compact \cite{P}.
Hence $G$ takes values in $L(E,E')$ \cite{D}.  Since $E$ has Property
(w), $G(A)$ is weakly compact for all $A \in \Sigma$.  If $\|G\|$ is
not continuous at $\emptyset$, then there are $(x''_n) \subseteq
U_{E''}$, a pairwise disjoint sequence $(A_n)$ in $\Sigma$, and
$\epsilon > 0$ such that $|G_{x''_n}(A_n)| > \epsilon$ for all $n$.
Choose $\Sigma$-measurable $U_E$-valued simple functions $(f_n)$ such
that supp\,$f \subseteq A_n$ and $\int\! f_ndG_{x''_n} > \epsilon$ for all $n$.
Since $U_E$ is $\sigma(E'',E')$-dense in $U_{E''}$ and $f_n$ is
simple, there exist $(x_n) \subseteq U_E$ such that $\int\! f_ndG_{x_n} >
\epsilon$.  Note that for all $x \in E$,
\[ \sum_n \lv \int\! f_ndG_x\rv \leq \sum_n|G_x|(A_n) \leq |G_x|(K) =
\|T'y\|. \]
Hence $S : E \to \ell^1$ defined by 
\[ Sx = \biggl(\int\! f_ndG_n\biggr)_n\]
is bounded.  Since $(Sx_n)(n) > \epsilon$ for all $n$, $(Sx_n)$ is not
weakly compact in $\ell^1$.  Using Proposition 2.a.2 of \cite{LT}, we
find a subsequence $(x_{n_m})$ such that $(Sx_{n_m})$ is equivalent to
the $\ell^1$ basis, and $[Sx_{n_m}]$ is complemented in $\ell^1$.
>From this it follows readily that $E$ contains a complemented copy of
$\ell^1$ as well.  This contradicts the fact that $E$ has Property
(w).  Since $G$ satisfies conditions (a)--(c) of Proposition \ref{op},
and $(K,E)$ is simple, the claim follows.\\

The proof that $C(K,E)$ has Property (w) proceeds analogously.
Let $T$ be an operator from  $C(K,E)$ into $M(K,E')$ represented by
the measure $G$, define
$T_x \in L(C(K),M(K,E''))$ as above by $T_xf = T(f\otimes x)$ for
all $x \in E, f \in C(K)$.  If $M(K,E')$ contains a copy of $c_0$,
then $\ell^1$ embeds complementably in $E$ \cite{SS2}, a
contradiction.  Hence $T_x$ does not fix a copy of $c_0$, and hence is
weakly compact \cite{P}.  Thus $G$ is $L(E,M(K,E'))$-valued \cite{D}.  For all
$A \in \Sigma$, let 
\[ S = G(A)'|_{C(K,E)} : C(K,E) \to E' . \]
$OBS$ is weakly compact by the claim above.  It is easy to see
that $S'x = G(A)x$ for all $x \in E$.  Thus $G(A) = S'|_{E}$ is weakly
compact.  Finally, if $\|G\|$ is not continuous at $\emptyset$, we
obtain, as in the proof of the claim, a pairwise disjoint sequence
$(A_n)$ in $\Sigma$,  a $\Sigma$-measurable
$U_E$-valued sequence of functions $(f_n)$ on $K$, a sequence $(g_n)
\subseteq U_{C(K,E)}$, and $\epsilon > 0$ such that supp$f_n \subseteq
A_n$ and $\int\! f_ndG_{g_n} > \epsilon$ for all $n$.  Now 
\[S : C(K,E)
\to \ell^1, \hspace{1em} Sg = \biggl(\int\! f_ndG_g\biggr)_n\]
is bounded.  Since $(Sg_n)(n) >
\epsilon$ for all $n$, we obtain as before that $\ell^1$ embeds
complementably into $C(K,E)$.  By \cite{SS2}, this implies that
$\ell^1$ embeds complementably into $E$, a contradiction.  Since the
pair $(K,E)$ is simple, we conclude that $T$ is weakly
compact.
\end{pf}
 
The pair $(K,E)$ is simple in either one of the following
situations:\\

\noindent(a) $E'$ and $E''$ both have the Radon-Nikod\'{y}m Property
\cite{BL};\\ 
\noindent(b) $K$ is a scattered compact \cite{C}.\\

\begin{cor}
If the pair $(K,E)$ satisfies one of the conditions {\rm (a)} or {\rm
(b)} listed
above, then $E$ has Property {\rm (w)} if and only if $C(K,E)$ does.
\end{cor}

\flushleft
\vspace{.5in}
Department of Mathematics\\National University of Singapore\\
Singapore 0511

\end{document}